\documentclass[11pt]{amsart}
 
\usepackage{amsfonts}
\usepackage{amssymb}
\usepackage{amscd}
 
\newcommand{\marginlabel}[1]%
  {\mbox{}\marginpar{\raggedleft\hspace{0pt}\bfseries\sf#1}}

\def\PP{{\textbf P}}

\def\KK{\mathcal{K}}
\def\cM{\mathcal{M}}

\def\Pic0{{\rm Pic}^0(X)}

\def\Mg{\overline{\mathcal{M}}_g}
\def\K{K_{\Mg}}
\def\K10{\overline{\mathcal{K}}}
\def\M10{\overline{\mathcal{M}}_{10}}

\theoremstyle{plain}

\newtheorem{theorem}{Theorem}[section]
\newtheorem{proposition}[theorem]{Proposition}

\newtheorem{corollary}[theorem]{Corollary}
\newtheorem{lemma}[theorem]{Lemma}

\theoremstyle{definition}

\newtheorem{remark}[theorem]{Remark}
\newtheorem{example}[theorem]{Example}
   
\newtheorem{conjecture}[theorem]{Conjecture}
\newtheorem{conjecture/question}[theorem]{Conjecture/Question}
\newtheorem{question}[theorem]{Question}  
\newtheorem{remark/definition}[theorem]{Remark/Definition} 
\newtheorem{terminology/notation}[theorem]{Terminology/Notation}
 
\theoremstyle{remark}

\begin{document}
 
\title{\bf Effective divisors on $\Mg$ and a counterexample to 
the Slope Conjecture}

\author[G. Farkas]{Gavril Farkas} 
\address{Department of Mathematics, University of Michigan,
525 East University, Ann Arbor, MI, 48109-1109}
\email{{\tt gfarkas@umich.edu}}
\thanks{Research of GF partially supported by the NSF Grant DMS-0140520}

\author[M. Popa]{Mihnea Popa}
\address{Department of Mathematics, Harvard University,
One Oxford Street, Cambridge, MA 02138}
\email{{\tt mpopa@math.harvard.edu}}
\thanks{Research of MP partially supported by the NSF Grant DMS-0200150}

\thanks{We thank Joe Harris and Rob Lazarsfeld for interesting discussions 
on this subject.}

\maketitle

\markboth{G. FARKAS and M. POPA}
{EFFECTIVE DIVISORS ON $\Mg$}

\section{\bf Introduction}

The purpose of this note is to prove two statements on the slopes 
of effective divisors on the moduli space of stable curves $\Mg$: 
first that the Harris-Morrison Slope Conjecture fails to 
hold on $\M10$ and second, that  in order 
to compute the slope of $\Mg$ for $g\leq 23$, one only has to look at the 
coefficients of the classes $\lambda$ and $\delta_0$ in the expansion
of the relevant divisors. 
The proofs are based on a general result providing inequalities
between the first few coefficients of effective divisors on $\Mg$.
We give the technical statements in what follows.

On $\Mg$ we denote by $\lambda$ the class of the Hodge line bundle, by
$\delta_0,\ldots, \delta_{[g/2]}$ 
the boundary divisor classes corresponding to singular stable curves
and by $\delta:=\delta_0+\cdots +\delta_{[g/2]}$ the total boundary.  
If $\mathbb E\subset \mbox{Pic}(\Mg)\otimes \mathbb R$ is the
effective cone, then following \cite{harris-morrison} we
define the \emph{slope function} 
$s:\mathbb E\rightarrow \mathbb R\cup \{\infty\}$ by the formula
$$s(D):= \mbox{inf }\{\frac{a}{b}:a,b>0 \mbox{ such that }
a\lambda-b\delta-D\equiv 
\sum_{i=0}^{[g/2]} c_i\delta_i,\mbox{ where }c_i\geq 0\}.$$
From the definition it follows that $s(D)=\infty$ unless $D\equiv
a\lambda-\sum_{i=0}^{[g/2]} b_i\delta_i$ 
with $a,b_i\geq 0$ for all $i$ (and it is well-known that $s(D)<\infty$ for any $D$
which is the closure of an effective divisor 
on $\mathcal{M}_g$). In the second case one has that 
$s(D)=a/{\mbox{min}}_{i=0}^{[g/2]}b_i.$
We denote by $s_g$ the slope of the moduli space 
$\Mg$, defined as $s_g:=\mbox{inf }\{s(D):D\in \mathbb E\}$.
The Slope Conjecture of Harris and Morrison predicts that 
$s_g\geq 6+12/(g+1)$ (cf. \cite{harris-morrison} Conjecture 0.1). 
This is known to hold for $g\leq 12, g\neq 10$ 
(cf. \cite{harris-morrison} and \cite{tan}).

Following \cite{cukierman-ulmer}, we consider the divisor $\KK$ on $\cM_{10}$ 
consisting of smooth curves lying on a $K3$ surface, and we denote by 
$\K10$ its closure in $\M10$. For any $g\geq 20$, we look at the locus in $\Mg$ 
of curves obtained by attaching a pointed curve of genus $g-10$ to a 
curve in $\K10$ with a marked point. This gives a divisor in $\Delta_{10}
\subset \Mg$, which we denote by $Z$.

The key point in what follows is that, based on the study of curves lying on $K3$ 
surfaces, one can establish inequalities involving a number of coefficients 
of any effective divisor coming from $\mathcal{M}_g$ in the expansion in 
terms of the generating classes.

\begin{theorem}\label{inequalities}
Let $D\equiv a\lambda-\sum_{i=0}^{[g/2]}b_i\delta_i$ be the closure in $\Mg$ of 
an effective divisor on $\mathcal{M}_g$.
\newline
\noindent (a) For $2\leq i\leq 9$ and $i=11$ we have $b_i \geq (6i+18)b_0 - (i+1)a$.
The same formula holds for $i=10$ if $D$ does not contain the divisor 
$Z\subset \Delta_{10}$.  
\newline
\noindent
(b) $(g\geq 20)$ If $D$ contains $Z$, then either $b_{10}\geq 78b_0-11a$ as above, or 
$b_{10}\geq (71.3866...)\cdot b_0 - (10.1980...)\cdot a$.
\newline
\noindent
(c) We always have that $b_1\geq 12b_0  - a$.
\end{theorem}

Part (a) (and (c)) of this theorem essentially only make more concrete results
--and the technique of intersecting with special Lefschetz pencils-- 
already existing in the literature mentioned above. Part (b) however is more 
involved: it requires 
pull-backs to $\overline{\mathcal{M}}_{10,1}$ and the intrinsic use of our 
partial knowledge about the divisor $\K10$, plus some facts about the 
Weierstrass divisor on $\overline{\mathcal{M}}_{g,1}$. Here we claim more originality. 

\begin{corollary}\label{bib0}
If $a/b_0\leq 71/10$, then $b_i \geq b_0$ for all $1\leq i\leq 9$. 
The same conclusion holds for $i=10$ if $a/b_0\leq 6.906\ldots$, and for $i=11$ if $a/b_0\leq 83/12$.
\end{corollary} 

Based on this we obtain that the divisor $\K10\subset \M10$ provides 
a counterexample to the Slope Conjecture.
Its class can be written as
$$\K10 \equiv a\lambda - b_0\delta_0  -\ldots - b_5\delta_5,$$
and by \cite{cukierman-ulmer} Proposition 3.5, we have $a=7$ and $b_0 = 1$. In view 
of Corollary \ref{bib0}, this information is sufficent to show that the 
slope of $\K10$ is smaller than the one expected based on the Slope Conjecture. 

\begin{corollary}\label{counterexample}
The slope of $\K10$ is equal to $a/b_0=7$, so strictly smaller than 
the bound $78/11$ predicted by the Slope Conjecture. In particular $s_{10}=7$ (since 
by \cite{tan} $s_{10}\geq 7$).
\end{corollary}

Theorem \ref{inequalities} also allows us to formulate (at least up to genus 
$23$, and conjecturally beyond that) the 
following principle: \emph{the slope $s_g$ of $\Mg$ is 
computed by the quotient $a/b_0$ of the relevant divisors}. We have 
more generally:

\begin{theorem}\label{slope}
For any $g\leq 23$, there exists  $\epsilon_g>0$ such that 
for any effective divisor $D$ on $\Mg$ with $s_g\leq s(D)\leq s_g + \epsilon_g$ we have 
$s(D) = a/b_0$, i.e. $b_0\leq b_i$ for all $i\geq 1$. 
\end{theorem}

\begin{conjecture}\label{slconj}
\emph{The  statement of the theorem holds in arbitrary genus.}
\end{conjecture}

\noindent
Theorem \ref{inequalities} is proved in \S2. Theorem \ref{slope} is proved
in \S3, where we also remark that the methods of the present paper give a
very quick proof of the fact that the Kodaira dimension of the universal curve 
$\overline{\mathcal{M}}_{g,1}$ is $-\infty$ for $g\leq 15$, $g\neq 13, 14$.
\noindent

\section{\bf Inequalities between coefficients of divisors}

Let $\mathcal{F}_g$ be the moduli space of canonically polarized 
$K3$ surfaces $(S,H)$ of genus $g$. We consider the $\PP^g$-bundle 
$\mathcal{P}_g=\{(S,C):C\in |H|\}$ over $\mathcal{F}_g$ which 
comes with a natural rational map 
$\phi_g:\mathcal{P}_g - ->\mathcal{M}_g$. By Mukai's results \cite{mukai1} 
and \cite{mukai2}, this
map is dominant if and only if $2\leq g\leq 9$ or $g=11$.
In this range $\mathcal{M}_g$ can be covered by curves corresponding to 
Lefschetz pencils of curves on $K3$ surfaces (cf. \cite{tan}). 
This is not true any more when $g=10$: in this case $\mbox{Im}(\phi_{10})$ 
is a divisor $\mathcal{K}$ in $\mathcal{M}_{10}$ (cf.
\cite{cukierman-ulmer} Proposition 2.2).
 
Given $2\leq i\leq 11$,
consider as above a Lefschetz pencil of curves of genus $i$ lying on a
general $K3$ surface of degree $2i-2$ in $\PP^i$. This gives rise to a
curve $B$ in the moduli space $\overline{\mathcal{M}}_i$.
Note that any such Lefschetz pencil, considered as a family of curves over $\PP^1$, has at least 
one section, since its base locus is nonempty.

\begin{lemma}\label{nos1}
We have the formulas $B\cdot \lambda =i+1$, $B\cdot \delta_0=
6i+18$ and $B\cdot \delta_j=0$ for $j\neq 0$.
\end{lemma}
\begin{proof}
The first two numbers are computed e.g. in \cite{cukierman-ulmer} 
Proposition
3.1, based on the formulas in \cite{griffiths-harris} pp. 508--509.
The last assertion is obvious since there are no reducible 
curves in a Lefschetz pencil.
\end{proof}

For each $g\geq i+1$, starting with the pencil $B$ in 
$\overline{\mathcal{M}}_i$ we can  construct a new
pencil $B_i$ in $\Mg$ in the following way: we fix a general pointed curve $(C,p)$ 
genus $g-i$. We then glue the curves in the 
pencil $B$ with $C$ at $p$, along one of the sections corresponding 
to the base points of the pencil.
We have that all such $B_i$ fill up $\Delta_i \subset \Mg$
for $i\neq 10$, and the divisor $Z\subset \Delta_{10}$ when $i=10$.

\begin{lemma}\label{nos2}
We have $B_i\cdot \lambda = i+1$, $B_i\cdot \delta_0 = 6i+18$, 
$B_i\cdot \delta_i = -1$ and $B_i\cdot \delta_j = 0$ for $j\neq 0, i$. 
\end{lemma}
\begin{proof}
This follows immediately from Lemma \ref{nos1} and from general principles, 
as explained in \cite{chang-ran} pp.271.
\end{proof}

\begin{proof}(\emph{of Theorem \ref{inequalities} (a), (c)})
(a) Let us fix $2\leq i\leq 11, i\neq 10$. Since $D$ is the closure of a
divisor coming from $\mathcal{M}_g$, it cannot contain the whole boundary $\Delta_i$. 
Thus we must have a pencil $B_i$ as above such that 
$B_i\cdot D\geq 0$. The same thing holds true for $i=10$ if we know
that $Z$ is not contained in $D$.
But by Lemma \ref{nos2} this is precisely the statement of this part.

\noindent
(c) This is undoubtedly well known. 
One follows the same procedure, but this time producing a pencil $B_1$ in 
$\Delta_1\subset \Mg$ by gluing a fixed pointed curve $(C,p)$ of genus $g-1$ to
a generic pencil of plane cubics along one of its $9$ sections. We have 
the well-known relations:
$$B_1\cdot\lambda =1, B_1\cdot \delta_0=12, B_1\cdot \delta_1 =-1 
{\rm ~and~} 
B_1\cdot \delta_j =0 {\rm~for~}j\neq 0,1.$$
The conclusion follows similarly, since we can find a $B_1$ such that 
$B_1\cdot D\geq0$.
\end{proof}

The study of the coefficient $b_{10}$ is more involved, since in $\mathcal{M}_{10}$ 
the Lefschetz pencils of curves on $K3$ surfaces only fill up a divisor.
We need some preliminaries on divisors on the universal curve 
$\overline{\mathcal{M}}_{g,1}$. Let  $\pi:\overline{\mathcal{M}}_{g,1}\rightarrow \Mg$
be the forgetful morphism.  The generators of $\mbox{Pic}(
\overline{\mathcal{M}}_{g,1})\otimes \mathbb Q$ are 
the tautological class  $\psi=c_1(\omega_{\pi})$, the boundary $\delta_0$, 
the Hodge class $\lambda$, 
and for $1\leq i\leq g-1$ the class $\delta_i$ corresponding to the 
locus of pointed curves
consisting of two components of genus $i$ and $g-i$ respectively 
with the marked point being on the genus $i$ component.

\begin{lemma}\label{push}(cf. \cite{arbarello-cornalba} \S1)
One has the following relations:
$$\pi_*(\lambda^2)=\pi_*(\lambda \cdot \delta_i)=\pi_*(\delta_{0}\cdot \delta_i)=0
{\rm~ for ~all~ } i=0,\ldots, g-1,\mbox{ }
\pi_*(\psi^2)=12\lambda-\delta,$$ 
$$\pi_*(\lambda\cdot \psi)=(2g-2)\lambda,\mbox{ } \pi_*(\psi\cdot \delta_0)=(2g-2)\delta_0,
\mbox{ } \pi_*(\psi\cdot \delta_i)=(2i-1)\delta_i {\rm~ for~} i\geq 1,$$
$$\pi_*(\delta_i^2)=-\delta_i {\rm~for~} 1\leq i\leq g-1,\mbox{ }\pi_*(\delta_i\cdot \delta_{g-i})=\delta_i,
{\rm~for~ } 1\leq i<g/2, {\rm ~and~ } $$
$$\pi_*(\delta_i\cdot \delta_j)=0 
{\rm~ for ~ all ~} i,j\geq 0 {\rm~ with ~} i\neq j, g-j.$$
\end{lemma}

\noindent
We consider the \emph{Weierstrass divisor} in $\mathcal{M}_{g,1}$
$$\mathcal{W}:=\{[C,p]\in \mathcal{M}_{g,1}:p\in C \mbox{ is a 
Weierstrass point}\},$$
and denote by $\overline{\mathcal{W}}$ its closure in $\overline{\mathcal{M}}_{g,1}$. 
Its class has been computed by Cukierman \cite{cukierman}:
$$\overline{\mathcal{W}}\equiv -\lambda+\frac{g(g+1)}{2}\psi-\sum_{i=1}^{g-1}
{g-i+1 \choose 2}\delta_i.$$

\begin{proposition}
If $\pi:\overline{\mathcal{M}}_{g,1}\rightarrow \overline{\mathcal{M}}_g$ 
is the forgetful morphism, then $\pi_*(\overline{\mathcal{W}}^2)$ is an effective 
divisor class on $\overline{\mathcal{M}}_g$.
\end{proposition}

\begin{proof} From the previous Lemma we have that $$\pi_*(\overline{\mathcal{W}}^2)
\equiv a\lambda-\sum_{i=0}^{[g/2]}b_i\delta_i,$$
where $a=g(g+1)(3g^2+g+2), b_0=g^2(g+1)^2/4$ while for $1\leq i<g/2$ we have $b_i=i(g-i)(g^3+3g^2+g-1)$. When $g$ is even $b_{g/2}=(8g^5+28g^3+33g^4+4g^2)/64.$
On the other hand we have expressions for the classes of distinguished 
geometric divisors on $\overline{\mathcal{M}}_g$: when $g+1$ 
is composite, by looking at Brill-Noether divisors one sees 
that the class
$$(g+3)\lambda-\frac{g+1}{6}\delta_0-\sum_{i=1}^{[g/2]}i(g-i)\delta_i$$
is effective (cf. \cite{eisenbud-harris} Theorem 1). When $g+1$ is 
prime one has to use the class of 
the Petri divisor, which gives a slightly worse estimate (cf. \cite{eisenbud-harris} Theorem 2).
In either case, by comparing the coefficients $a,b_i$ above with those of 
these explicit effective classes, one obtains an 
effective representative for $\pi_*(\overline{\mathcal{W}}^2)$.
For instance when $g+1$ is composite it is enough to check that $b_0/a\leq (g+1)/(6g+18)$ 
and that $b_i/a\leq i(g-i)/(g+3)$ for $i=1,\ldots,[g/2]$, which is immediate.
\end{proof}

\begin{corollary}\label{wei}
Let $D$ be any effective divisor class on $\overline{\mathcal{M}}_{g,1}$. 
Then $\pi_*(\overline{\mathcal{W}}\cdot D)$ is an effective class on 
$\overline{\mathcal{M}}_g$. 
\end{corollary}

\begin{proof}(\emph{of Theorem \ref{inequalities} (b)})
Assume that $b_{10}<78b_0-11a$. 
We consider the map 
$$j:\overline{\mathcal{M}}_{10,1}\longrightarrow \Mg$$  
obtained by attaching a fixed general pointed 
curve of genus $g-10$ to any curve of genus $10$ with a marked point. 
Our assumption says that $R \cdot j^*(D)<0$,
where $R\subset \overline{\mathcal{M}}_{10,1}$ denotes the
curve in the moduli space coming from a Lefschetz pencil of 
pointed curves of genus $10$ on a general $K3$ surface. 
We can write $j^*(D)=m\pi^*(\overline{\mathcal{K}})+E,$ where $E$ is an effective 
divisor not containing $\pi^*(\overline{\mathcal{K}})$ and 
$m\in \mathbb Z$ is such that
\begin{equation}\label{mm}
m\geq -R\cdot j^*(D) =-11a+78b_0-b_{10}>0.
\end{equation}
Note that we have the standard formulas 
$j^*(\lambda)=\lambda, j^*(\delta_0)=\delta_0$ and $j^*(\delta_{10})=-\psi$
(cf. \cite{arbarello-cornalba}), while $\overline{\mathcal{K}}\equiv 7\lambda-\delta_0-\cdots,$ 
hence
$$E\equiv (a-7m)\lambda+b_{10}\psi+(m-b_0)\delta_0+
{\rm~(other~ boundaries)}$$ 
is an effective class on 
$\overline{\mathcal{M}}_{10,1}$. By applying Corollary \ref{wei}
it follows that $\pi_*(\overline{\mathcal{W}}\cdot E)$ is an 
effective class on $\overline{\mathcal{M}}_{10}$. An easy calculation 
using 
Lemma \ref{push} shows that 
$$\pi_*(\overline{\mathcal{W}}\cdot E)\equiv \bigl(642b_{10}+990(a-7m)\bigr)
\lambda-55\bigl(b_{10}
+18(b_0-m)\bigr)\delta_0-\cdots.$$
We now use that for every effective divisor on $\Mg$ the 
coefficient $a$ of $\lambda$ is nonnegative \footnote{For the reader's 
convenience we recall that this follows immediately from the fact that 
$B\cdot \lambda>0$ for any curve $B\subset \Mg$ such that $B\cap
\mathcal{M}_g\neq \emptyset$, 
while there is always a complete curve in
$\mathcal{M}_g$ passing through a general point.}.
From the previous formula we get an inequality which combined with 
(\ref{mm}) yields, after a simple computation
$$b_{10}\geq (71.3866...)\cdot b_0 - (10.1980...)\cdot a.$$
\end{proof}

\begin{question}
Do we always have the inequality $b_{10}\geq 78b_0-11a$?
\end{question}

We conclude with some examples where these inequalities can be 
checked directly and are sometimes sharp.

\begin{example}
When $g+1$ is composite, if $r, d>0$ are such that $g+1=(r+1)(g-d+r)$, then 
the locus of curves of genus $g$ carrying a $\mathfrak g^r_d$  is a divisor 
with class (cf. \cite{eisenbud-harris}, Theorem 1): 
$$\overline{\mathcal{M}}_{g,r}^d \equiv c\big((g+3)\lambda -\frac{g+1}{6} \delta_0
-\sum_{i=1}^{[g/2]}i(g-i)\delta_i\big),$$
where $c$ is a  constant depending on $g,r$ and $d$. 
A simple calculation shows that the inequalities in Theorem 
\ref{inequalities}
are satisfied. Moreover, they are sharp for $i=1$ and $i=2$, 
for any genus $g$. 
\end{example}

\begin{example}
A similar behavior is exhibited by the divisor 
of curves on $K3$ surfaces $\K10\subset \M10$, where we have $b_1=5$ 
and $b_2=9$ which again gives equality in the first 
two inequalities in Theorem \ref{inequalities}. Note that this 
follows by the method of \cite{eisenbud-harris} \S2 if we show that 
the pull-back of $\K10$ to $\overline{\mathcal{M}}_{2,1}$ is 
supported on the Weierstrass divisor. We will obtain this in 
the forthcoming paper \cite{future}, as a special case of a more general 
study of degenerations of spaces of sections of rank two vector bundles on curves.
The same study will show a substantial difference between the geometry 
of $\K10$ and that of the Brill-Noether divisors, 
namely that the image of the natural \lq \lq flag" map from $\overline{\mathcal{M}}_{0,g}$
to $\Mg$ is contained in the $K3$-locus for any $g$. Thus one cannot 
use the method of \cite{eisenbud-harris} \S3 in order to determine more 
coefficents of $\K10$.
\end{example}

\section{\bf Slopes of divisors and further remarks}

The inequalities established in Theorem \ref{inequalities} allow 
us to show that, at least up to genus $23$, if the slope of an
effective divisor is sufficiently small, then it is computed by the 
ratio $a/b_0$.

\begin{proof}(\emph{of Theorem \ref{slope}})
When $g$ is such that $g+1$ is composite, we have that $s_g\leq
6+12/(g+1)$ 
(this being the slope of any Brill-Noether divisor). When $g$ is even, one has the estimate
$s_g\leq \frac{2(3g^2+13g+2)}{g(g+2)}$ (this being the slope of the
Petri divisor, 
cf. \cite{eisenbud-harris} Theorem 2).
It follows that  for any $g\leq 23$ there exists
a positive number $\epsilon_g$ such that 
$$s_g+\epsilon_g\leq 6+ \frac{11}{i+1} {\rm~for~all~} i\leq [g/2](\leq 11).$$

Assume first that $2\leq i\leq 9$ or $i=11$. Then by Theorem 
\ref{inequalities}(a)
we know that $b_i \geq (6i+18)b_0 - (i+1)a$, and so certainly $b_i\geq b_0$ 
if $s(D)\leq 6+ \frac{11}{i+1}$. 
For $i=10$ we apply \ref{inequalities}(b): if the inequality $b_{10}\geq
78b_0 -11a$ holds, then the argument is identical. If not, we have the 
inequality $b_{10}\geq (71.3866...)\cdot b_0 - (10.1980...)\cdot a$. 
Thus $b_{10}\geq b_0$
as soon as the inequality $a/b_0\leq 6.9$ is satisfied. But for 
$g\geq 20$ 
the inequality $s_g < 6.9$ holds, based on the same estimates as above.
For $i=1$, the condition 
is even weaker because of the formula $b_1\geq 12b_0-a$ in \ref{inequalities}(c). 
Thus the slope of $D$ is computed by $a/b_0$.
\end{proof}

\begin{remark}
Let us consider again the curve $B\subset \Mg$ corresponding to a Lefschetz
pencil of curves of genus $g$ on a general $K3$ surface (cf. Lemma \ref{nos1}).
Since $B\cdot \delta_0/B\cdot \lambda=6+12/(g+1)$ which is the conjectured
value of $s_g$, it follows that the nefness  
\footnote{Recall that, slightly abusively, a curve $B$ on a projective 
variety $X$ is called \emph{nef} 
if $B\cdot D\geq 0$ for every efective Cartier divisor $D$ on $X$.} 
of $B$ would be a sufficient condition 
for the Slope Conjecture to hold in genus $g$. Moreover, Theorem 
\ref{slope} and Corollary \ref{bib0} imply 
that for $g\leq 23$ the Slope Conjecture in genus $g$ is
equivalent to $B$ being a nef curve. The conjecture fails for $g=10$ because
$B\cdot \overline{\mathcal{K}}=-1$.
\end{remark}

\begin{remark} 
An amusing consequence of Proposition \ref{wei} is that 
the Kodaira dimension of the universal curve $\mathcal{M}_{g,1}$ 
is $-\infty$ for all $g\leq 15$, with $g\neq 13,14$ ( of course this 
can be proved directly when $g\leq 11$).
Indeed, if we assume that the canonical class $K_{\overline{\mathcal{M}}_{g,1}}\equiv 
13\lambda+\psi-3(\delta_1+\delta_{g-1})-2\sum_{i=2}^{g-2}\delta_i$  is effective on 
$\overline{\mathcal{M}}_{g,1}$, then by Proposition \ref{wei}
the class $D:=\pi_*(K_{\overline{\mathcal{M}}_{g,1}}\cdot \overline{\mathcal{W}})$ 
is effective on $\overline{\mathcal{M}}_g$. 
It turns out that
$s(D)=\frac{2(13g^3+6g^2-9g+2)}{g(g+1)(4g+3)}$, and from the 
definition of the slope of $\mathcal{M}_g$ we have 
that $s(D)\geq s_g$. But this contradicts  the estimates on $s_g$
from \cite{tan} and \cite{chang-ran}.
\end{remark}

\end{document}